\documentclass[12pt]{article}

\usepackage{latexsym}
\usepackage{amsfonts}
\usepackage{amsmath}

 \newcommand{\IR}{{\Bbb R}}

\newcommand{\nn}{{\nonumber}}

\newtheorem{theorem}{Theorem}[section]

\begin{document}

\title{Central limit theorems for a class of irreducible multicolor urn models}

\author{Gopal K. Basak\thanks{
Department of Mathematics, University of Bristol, University Walk, Bristol, BS8 1TW, UK,  and Stat-Math Unit, Indian Statistical Institute, 203 B. T. Road,
Kolkata 700108, INDIA, Email Address:
magkb@bris.ac.uk and gkb@isical.ac.in}
 \and Amites Dasgupta\thanks{
Stat-Math Unit, Indian Statistical Institute, 203 B. T. Road,
Kolkata 700108, INDIA, Email Address: amites@isical.ac.in} \\
}
\date{}

\maketitle

\begin{abstract}

We take a unified approach to central limit theorems for a class of irreducible urn models with constant replacement matrix. Depending on the eigenvalue, we consider appropriate linear combinations of the number of balls of different colors. Then under appropriate norming the multivariate distribution of the weak limits of these linear combinations is obtained and independence and dependence issues are investigated.  

\end{abstract}

{\bf Keywords and Phrases :} Central limit theorem, Urn models, Martingale.

{\bf 2000 Subject classification:} Primary: 60F17; Secondary: 60J30, 60G15, 60G45. 

\newpage

\section{Introduction}

Consider a four-color urn model in which the replacement matrix is actually a stochastic matrix $\mathbf{R}$ in the manner of Gouet \cite{gouet2}. That is, we start with one ball of any color, which is the $0$-th trial. Let $\mathbf{W}_n$ denote the column vector of the number of balls of the four colors upto the $n$-th trial, where the components of $\mathbf{W}_n$ are nonnegative real numbers.  Then a color is observed by random sampling from a multinomial distribution with probabilities $(1/(n+1)) \mathbf{W}_n$. Depending on the color that is observed, the corresponding row of $\mathbf{R}$ is added to $\mathbf{W}_n^\prime$ and this gives $\mathbf{W}_{n+1}^\prime$. A special case of the main theorem of Gouet \cite{gouet2} is that if the stochastic matrix $\mathbf{R}$ is irreducible, then $(1/(n+1)) \mathbf{W}_n^\prime$  converges a.s. to the stationary distribution $\mathbf{\pi}$ of the irreducible stochastic matrix $\mathbf{R}$ (it should be carefully noted that the multicolor urn model is vastly different from the Markov chain evolving according to the transition matrix equal to the stochastic matrix $\mathbf{R}$). Suppose the nonprincipal eigenvalues of $\mathbf{R}$ satisfy $\lambda_1 < 1/2, \lambda_2 = 1/2, \lambda_3 > 1/2$ respectively, which are assumed to be real (and hence lie in $(-1, 1)$), and $\xi_1, \xi_2, \xi_3$ be the corresponding eigenvectors. Using $\mathbf{\pi} \xi_i = \mathbf{\pi} \mathbf{R} \xi_i = \lambda_i \mathbf{\pi} \xi_i$ it is seen that $(1/(n+1)) \mathbf{W}_n^\prime \xi_i \rightarrow 0$. Thus central limit theorems are the next interesting statistical results. 

In this article we consider the joint limiting distribution of $(X_n, Y_n, Z_n)$ where 
\begin{equation}
X_n = \frac{\mathbf{W}_n^\prime \xi_1}{\sqrt{n}}, 
Y_n = \frac{\mathbf{W}_n^\prime \xi_2}{\sqrt{n \log n}},
Z_n = \frac{\mathbf{W}_n^\prime \xi_3}{\Pi_0^{n-1} (1 + \frac{\lambda_3}{j+1})} . \label{prai} 
\end{equation} 
Special cases of this result are known from Freedman \cite{freedman}, Gouet \cite{gouet} and Smythe \cite{smythe} and Bai and Hu \cite{bai1}. Freedman \cite{freedman}, as well as Gouet \cite{gouet}, consider two color urn, so that there is only one eigenvector and the corresponding nonprincipal eigenvalue can be one of the three types. The identification of the norming rates is due to Freedman \cite{freedman}, his technique being the method of moments. Smythe \cite{smythe} considers multicolor urn, but all the nonprincipal eigenvalues (or their real parts) are assumed to be $< 1/2$. Recently Bai and Hu \cite{bai1} have considered the case when all the nonprincipal eigenvalues (or their real parts) are $\leq 1/2$. However to our knowledge, the joint limit when eigenvalues of all the three types occur has not been considered. The interesting feature of this case, which will be clear from the proof, is the differences in the behaviors of the differences of the three components. 

For the above four color set up our main result is 
\begin{theorem} \label{th.1}
$(X_n, Y_n, Z_n)$ converges in distribution to $(X, Y, Z)$ where \newline  $X, Y, Z$ are independent, $X$ and $Y$ are (independent) normals with zero means. The convergence of $Z_n$ to $Z$ is also  in the almost sure sense.
\end{theorem}
The variances of $X$ and $Y$ are identified in the proof. However our proof says nothing about the distribution of $Z$, apart from $E Z = 0$. Some features of this $Z$ in a two color case are discussed in Freedman \cite{freedman}. We also need to point out the connection of theorem 1.1 with the available results in the literature. The available results actually consider norming the vector $(\mathbf{W}_n - E \mathbf{W}_n)$ and not the linear combinations from the eigenvectors. Now the eigenvectors $\xi_1, \xi_1, \xi_3$ and the principal eigenvector $\mathbf{u} = (1, 1, 1, 1)^\prime$ span $\IR^4$, so that any linear combination can be expressed in terms of them. But $\mathbf{W}_n^\prime \mathbf{u} = n + 1$, so its effect cancels out after the expectation is subtracted and we are left with the linear combinations corresponding to $\xi_1, \xi_2, \xi_3$. The available results in the literature divide  $(\mathbf{W}_n - E \mathbf{W}_n)$ by the largest rate, and in the case the real part of the nonprincipal eigenvalues is less than or equal to 1/2 (actually the rate in that case may be different from $\sqrt{n \log n}$ as will be clear in the later sections) derive asymptotic normality, see e.g. Bai et al. \cite{bai1}. 

We have stated our theorem for the four color model for the sake of notational simplicity in the proof. The theorem also extends to situations (with more than four colors) where there are more than one eigenvalue(s) of any one or more of the three types. These extensions involve the same technique, but require more calculations related to the Jordan form of the replacement matrix. So we have sketched some of these extensions in separate sections. These sections discuss the main theorem in increasing generality along with development of suitable notation, and we have indicated the generalizations inside these sections. First, all the eigenvalues are considered to be real, the Jordan form thus involves only real vectors. Next, the eigenvalues can be complex, so the Jordan form involves complex vectors and we deal with the real and imaginary parts of these vectors. Another interesting feature of these later sections dealing with the Jordan form is the appearance of nilpotent and rotation matrices. The final result is given as theorem 5.1.

The proof of theorem 1.1 for the above four color set up is given in the next section. It employs an iteration technique involving conditional characteristic functions (an example of these iterations occurs in example 2, Pp. 79-80 of Basak et al. \cite{basak}). We have written this proof in detail, however the proofs for the generalizations of the main theorem are only sketched in later sections as the ideas are the same.      

\section{Proof of theorem 1.1}

Before starting the proof we collect a few computational details.
The column vector of the indicator functions of balls of different colors obtained from the $n+1$-st trial is denoted by $\mathbf{\chi}_{n+1}$. It is clear that $E\{ \mathbf{\chi}_{n+1}|{\cal F}_n \} = (1/(n+1)) \mathbf{W}_n$, where ${\cal F}_n$ denotes the $\sigma$-field of observations upto the $n$-th trial. This notation leads to \begin{equation}
\mathbf{W}_{n+1}^\prime \xi_i = \mathbf{W}_n^\prime \xi_i + \mathbf{\chi}_{n+1} \mathbf{R} \xi_i = \mathbf{W}_n^\prime \xi_i + \lambda_i \mathbf{\chi}_{n+1} \xi_i. \label{cru}
\end{equation}
For the purpose of iteration we shall use a decomposition of the components of the Markov chain $(X_{n+1}, Y_{n+1}, Z_{n+1})$ illustrated with the first component as follows: \[ X_{n+1} =  E\{ X_{n+1}|{\cal F}_n\} + (X_{n+1} - E\{X_{n+1} |{\cal F}_n\}).\] The first term will be expressed in terms of $X_n$ and the second term is the martingale difference that will play an important role in our proof in analogy with the calculations for the central limit theorem for i.i.d. random variables. 

To write the first term in terms of $X_n$ ($Y_n, Z_n$ respectively) we shall use the following approximations 
\begin{eqnarray*}
(1 + 1/n)^{-1/2} &=& 1 - \frac{1}{2n} + O(\frac{1}{n^2}), \\
\frac{\log n}{\log (n+1)} &=& \frac{\log n}{\log n + 1/n + O(1/n^2)} \\
&=& \frac{1}{1 + (1/n\log n) + O(1/n^2 \log n)}, \\
\sqrt{\frac{n \log n}{(n+1)\log (n+1)}} &=& \{ 1 - \frac{1}{2n} + O(\frac{1}{n^2})\} \{ 1 - \frac{1}{2n \log n} + O(\frac{1}{n^2}) \} \\
&=& 1 - \frac{1}{2n} - \frac{1}{2n \log n} + O(\frac{1}{n^2}), \\
\Pi_0^{n - 1} (1 + \lambda_3/(j+1)) &\sim& n^{\lambda_3}.
\end{eqnarray*}
 Using these and the conditional expectation of (\ref{cru}) it follows that:
\begin{eqnarray}
E \{ X_{n+1}|{\cal F}_n\} &=& X_n (1 - \frac{1/2 - \lambda_1}{n}) + X_n O(1/n^2), \nonumber \\
E\{ Y_{n+1}|{\cal F}_n\} &=& Y_n ( 1  - \frac{1}{2n \log n}) + Y_n O(1/n^2), \nonumber \\ 
E\{ Z_{n+1}|{\cal F}_n\} &=& Z_n, \label{cial}
\end{eqnarray}
the second of which crucially uses $\lambda_2 = 1/2$.
 Now let us look at the martingale difference terms which are:
\begin{eqnarray}
M_{1,n+1} = X_{n+1} - E\{X_{n+1} |{\cal F}_n\} &=& \lambda_1 \frac{\chi_{n+1}^\prime \xi_1}{\sqrt{n+1}} - \frac{\lambda_1}{n+1} \sqrt{\frac{n}{n+1}} X_n,  \nonumber \\
M_{2,n+1} = Y_{n+1} - E\{Y_{n+1} |{\cal F}_n\} &=& \lambda_2 \frac{\chi_{n+1}^\prime \xi_2}{\sqrt{(n+1) \log (n+1)}} \nonumber \\
&&  - \frac{\lambda_2}{n+1} Y_n \sqrt{\frac{n \log n}{(n+1) \log (n+1)}} \nonumber \\
M_{3, n+1} = Z_{n+1} - E\{Z_{n+1} |{\cal F}_n\} &=&  \lambda_3 \frac{\chi_{n+1}^\prime \xi_3}{\Pi_0^n(1 + \frac{\lambda_3}{j+1})} -\frac{\frac{\lambda_3}{n+1}}{1 + \frac{\lambda_3}{n+1}} Z_n. \label{prime}
\end{eqnarray}
 It will be seen that the part involving $\chi_{n+1}^\prime \xi_i$ plays the significant role in the second moment calculations.

\subsection{Main idea of the proof}

Now we are ready to start the proof of theorem 1.1. 

{\bf Step A:} Using (\ref{cial}) and the inequality
$|e^{ix} - 1| \leq 2|x|$ for real number $x$, and remembering that $|\mathbf{W}_n^\prime \xi_i| \leq c n$, so that $X_n/\sqrt{n}, Y_n/\sqrt{n}, Z_n/n^{1 - \lambda_3}$ are bounded, we can expand $e^{i t_1 X_n O(1/n^2) + i t_2 Y_n O(1/n^2)}$ to get
\begin{eqnarray}
&& \Big{|} E \{ e^{i (t_1 X_{n+1} + t_2 Y_{n+1} + t_3 Z_{n+1})}|{\cal F}_n \} \nonumber \\
&& - e^{i \{ t_1 (1 - \frac{\frac{1}{2} - \lambda_1}{n}) X_{n} + t_2 (1 - \frac{1}{2n \log n})  Y_{n} + t_3 Z_{n}\}}  E \{ e^{i (t_1 M_{1,n+1} + t_2 M_{2,n+1} + t_3 M_{3,n+1})}|{\cal F}_n \} \Big{|} \nonumber \\
&& \leq 2 (|t_1| |X_n| + |t_2| |Y_n|) O(1/n^2) \nonumber \\
&& \leq const. \frac{1}{n^{3/2}}, \label{ke} 
\end{eqnarray}
for $n$ sufficiently large, say $n \geq n_0$. 

{\bf Step B:} Now we want to approximate $E \{ e^{i (t_1 M_{1,n+1} + t_2 M_{2,n+1} + t_3 M_{3,n+1})}|{\cal F}_n \}$ by 
\begin{equation}
 e^{- \frac{t_1^2}{2} \lambda_1^2 \frac{<\pi, \xi_1^2 >}{n+1} - \frac{t_2^2}{2}\lambda_2^2 \frac{<\pi, \xi_2^2 >}{(n+1)\log (n+1)} }.\label{polo}
\end{equation} 
We use the inequality $|e^{ix} - 1 - ix + \frac{1}{2} x^2| \leq const. |x|^3$ along with the observation that the martingale differences of (\ref{prime}) are bounded by $const. /\sqrt{n}$, $const./\sqrt{n \log n}$ and $const./n^{\lambda_3}$ respectively (we approximate $\Pi_0^n (1 + \lambda_3/(i+1)) \sim n^{\lambda_3}$). This gives 
\begin{eqnarray}
&& \Big{|} E \{ e^{i (t_1 M_{1,n+1} + t_2 M_{2,n+1} + t_3 M_{3,n+1})}|{\cal F}_n \} \nonumber \\
&& - (1 - \frac{1}{2} E \Big{\{} (t_1^2 M_{1, n+1}^2 +  t_2^2 M_{2, n+1}^2 + t_3^2 M_{3, n+1}^2 \nonumber \\
&& + t_1 t_2 M_{1, n+1} M_{2, n+1} + t_1 t_3 M_{1, n+1} M_{3, n+1} + t_2 t_3 M_{2, n+1} M_{3, n+1})|{\cal F}_n \Big{\}}) \Big{|} \nonumber \\
&& \leq const. \frac{1}{n^{3/2}}, \label{re} 
\end{eqnarray}
for $n \geq n_0$. 

To achieve (\ref{polo}) a detailed study of the terms of (\ref{re}) is necessary. We denote by $\xi_i \xi_j$ the vector whose components are products of the corresponding components of $\xi_i$ and $\xi_j$. Remembering that $\chi_{n+1}$ consists of indicator functions of observations of balls of different colors, we get
\begin{eqnarray}
E (M_{1, n+1}^2 |{\cal F}_n) &=& \lambda_1^2 \frac{<\pi, \xi_1^2 >}{n+1} \nonumber \\
&& + \Big{\{} \lambda_1^2 \frac{< \frac{\mathbf{W}_n^\prime}{n+1} - \pi, \xi_1^2 >}{n+1} -  \lambda_1^2 \frac{n}{(n+1)^3} X_n^2 \Big{\}}, \nonumber \\ 
E (M_{2, n+1}^2 |{\cal F}_n) &=& \lambda_2^2 \frac{<\pi, \xi_2^2 >}{(n+1)\log (n+1)} \nonumber \\
&& + \Big{\{} \lambda_2^2 \frac{< \frac{\mathbf{W}_n^\prime}{n+1} - \pi, \xi_2^2 >}{(n+1)\log (n+1)}  -  \lambda_2^2 \frac{n \log n}{(n+1)^3 \log (n+1)} Y_n^2 \Big{\}}, \nonumber \\
E (M_{3, n+1}^2 |{\cal F}_n) &=& \lambda_3^2 \frac{<\pi, \xi_3^2 >}{(\Pi_0^n(1 + \frac{\lambda_3}{j+1}))^2} \nonumber \\
&& + \Big{\{} \lambda_3^2 \frac{<  \frac{\mathbf{W}_n^\prime}{n+1} - \pi, \xi_3^2 >}{(\Pi_0^n(1 + \frac{\lambda_3}{j+1}))^{2}} -  \frac{\lambda_3^2}{(n+1)^2 (1 + \frac{\lambda_3}{n+1})^2} Z_n^2 \Big{\}}, \nonumber 
\end{eqnarray}
\begin{eqnarray}
E (M_{1, n+1} M_{2, n+1}|{\cal F}_n) &=& \lambda_1 \lambda_2 \frac{<\pi , \xi_1 \xi_2 >}{\sqrt{n+1} \sqrt{ (n+1) \log (n+1) }} \nonumber \\
&& + \Big{\{} \lambda_1 \lambda_2 \frac{<\frac{\mathbf{W}_n^\prime}{n+1} -  \pi , \xi_1 \xi_2 >}{\sqrt{n+1} \sqrt{ (n+1) \log (n+1) }} \nonumber \\
&&  - \lambda_1 \lambda_2 \frac{n \sqrt{\log n}}{(n+1)^3 \sqrt{\log (n+1)}} X_n  Y_n \Big{\}},\nonumber \\
E (M_{1, n+1} M_{3, n+1}|{\cal F}_n) &=& \lambda_1 \lambda_3 \frac{<\pi, \xi_1 \xi_3>}{\sqrt{n+1} (\Pi_0^n(1 + \frac{\lambda_3}{j+1}))} \nonumber \\
&& + \Big{\{} \lambda_1 \lambda_3 \frac{<\frac{\mathbf{W}_n^\prime}{n+1} -  \pi, \xi_1 \xi_3>}{\sqrt{n+1} (\Pi_0^n(1 + \frac{\lambda_3}{j+1}))} \nonumber \\
&& - \lambda_1 \lambda_3 \sqrt{\frac{n}{n+1}} \frac{\frac{\lambda_3}{n+1}}{(n+1) (1 + \frac{\lambda_3}{n+1})} X_n Z_n \Big{\}}, \nonumber \\
E (M_{2, n+1} M_{3, n+1}|{\cal F}_n) &=& \lambda_2 \lambda_3 \frac{<\pi, \xi_2 \xi_3>}{\sqrt{(n+1) \log (n+1) } (\Pi_0^n(1 + \frac{\lambda_3}{j+1}))} \nonumber \\
&& + \Big{\{} \lambda_2 \lambda_3 \frac{<\frac{\mathbf{W}_n^\prime}{n+1} - \pi, \xi_2 \xi_3>}{\sqrt{(n+1) \log (n+1) } (\Pi_0^n(1 + \frac{\lambda_3}{j+1}))} \nonumber \\
&& - \lambda_2 \lambda_3 \sqrt{\frac{n \log n}{(n+1) \log (n+1)}} \frac{\frac{\lambda_3}{n+1}}{(n+1) (1 + \frac{\lambda_3}{n+1})} Y_n Z_n \Big{\}}. \label{kere}
\end{eqnarray}
If $\sigma^2$ is small then we know that $|1 - \frac{\sigma^2}{2} - e^{-\sigma^2/2}| \leq const. \sigma^4$. Using this on the constant terms of the first two equations of (\ref{kere}) we get 
\begin{eqnarray}
&& \Big{|} 1 - \frac{1}{2} t_1^2 \lambda_1^2 \frac{<\pi, \xi_1^2 >}{n+1} - \frac{1}{2} t_2^2 \lambda_2^2 \frac{<\pi, \xi_2^2 >}{(n+1)\log (n+1)} \nonumber \\
&& - e^{- \frac{t_1^2}{2} \lambda_1^2 \frac{<\pi, \xi_1^2 >}{n+1} - \frac{t_2^2}{2}\lambda_2^2 \frac{<\pi, \xi_2^2 >}{(n+1)\log (n+1)} } \Big{|} \nonumber \\
&& \leq const. \frac{1}{(n+1)^2}. \label{vate} 
\end{eqnarray}

{\bf Step C:} Combining (\ref{ke}), (\ref{re}) and (\ref{vate}) we get the following basic inequality:
\begin{eqnarray}
&& \Big{|} E \{ e^{i (t_1 X_{n+1} + t_2 Y_{n+1} + t_3 Z_{n+1})}|{\cal F}_n \} \nonumber \\
&& - e^{i \{ t_1 (1 - \frac{\frac{1}{2} - \lambda_1}{n}) X_{n} + t_2 (1 - \frac{1}{2n \log n})  Y_{n} + t_3 Z_{n}\}} \nonumber \\
&& \times e^{- \frac{t_1^2}{2} \lambda_1^2 \frac{<\pi, \xi_1^2 >}{n+1} - \frac{t_2^2}{2}\lambda_2^2 \frac{<\pi, \xi_2^2 >}{(n+1)\log (n+1)} }\Big{|} \nonumber \\
&& \leq const. \frac{1}{n^{3/2}} + R_n, \label{chess} 
\end{eqnarray}
where we use $R_n$ to denote the sum of the other constant terms and random terms from the right of (\ref{kere}) which have not been used in (\ref{vate}) (this is also multiplied by exponentials of imaginary quantities, but those are bounded by 1 and will not make any difference). We also use the notation \[ C_n = - \frac{t_1^2}{2} \lambda_1^2 \frac{<\pi, \xi_1^2 >}{n+1} - \frac{t_2^2}{2}\lambda_2^2 \frac{<\pi, \xi_2^2 >}{(n+1)\log (n+1)}.\]
We then condition again on ${\cal F}_{n-1}$ and iterate backwards. While doing so, in the exponent the coefficients of $t_i$ change as above, we get a sum of $C_{n - j}$'s in the exponent, and following iteration of (\ref{chess}) on the right we get a sum of conditional expectations of $R_n$'s and $const. \sum_{n_0}^n 1/(j+1)^{3/2}$. Note that the iteration from $n+1$ to $n$ has changed the coefficient of $X_n$ and $Y_n$, and these are assumed to be incorporated in $C_{n-1}$ and $R_{n-1}$, and so on. $R_{n-j}$ also involves terms like $e^{C_{n - j+1} + \cdots + C_n}$, but it will be seen from steps 1 and 2 that these terms are bounded uniformly and will be  absorbed in the $const.$ term in (\ref{abar}). We should mention here that the constant term in (\ref{ke}), (\ref{re}) and (\ref{vate}) and finally (\ref{chess}) can be taken independently of this iteration because during the iteration the coefficients of $t_1$ and $t_2$ decrease. 

The main idea of the proof is to iterate the (conditional) characteristic function backwards upto a sufficiently large $n_0$, and first make $n \rightarrow \infty$. This will make the sum of $C_n$'s independent of $n_0$, and the sum of the conditional expectations of the $R_n$'s given ${\cal F}_{n_0}$ will be bounded by a random variable (which depends on the fixed $n_0$).  Taking expectation of the conditional characteristic function we get the characteristic function. Then we let $n_0 \rightarrow \infty$, and a further argument gives us the characteristic function. Before we do this we provide a few ingredients of the proof in a separate subsection.

\subsection{Important limits and estimates}

So {\em assume we have iterated backwards upto a sufficiently large} $n_0$. For ease of exposition we divide the calculations into a few steps. In step 1 we concentrate on the nonrandom terms corresponding to $t_1^2$ and $t_2^2$, which gives the form of the characteristic function corresponding to $X_n$ and $Y_n$. In step 2 we consider the other nonrandom terms, then in step 3 we handle the random (second bracketed) terms. Steps 2 and 3 contribute to the sum of $R_n$'s.

{\bf Step 1:} The calculations here will go into $C_n$. They come from the first (nonrandom) terms of the first two equations on the right of (\ref{kere}). Because of the presence of the term $(1 - \frac{\frac{1}{2} - \lambda_1}{n})$ in the characteristic function, it is seen that after iterating backwards upto $n_0$, the (nonrandom part of the) coefficient of $-(1/2)t_1^2$ is 
\[ \sum_{n_0}^n f_{n-j+1} \lambda_1^2 \frac{<\pi, \xi_1^2 >}{j+1} \]
where \[ f_{n-j+1} = \Pi_{i = j+1}^n (1 - \frac{\frac{1}{2} - \lambda_1}{i})^2 \]
As $n \rightarrow \infty$, the above sum clearly goes to 
\begin{equation}
\lambda_1^2 < \pi, \xi_1^2> \int_0^\infty e^{-(1 - 2\lambda_1)x} dx. \label{pharao}
\end{equation}
Similarly because of the presence of $(1 - \frac{1}{2n \log n})$ in the characteristic function, after iterating backwards upto $n_0$, the (nonrandom part of the) coefficient of $- \frac{1}{2} t_2^2$ is  \[ \sum_{n_0}^n g_{n-j+1} \lambda_2^2 \frac{<\pi, \xi_2^2 >}{(j+1) \log (j+1) } \]
where \[ g_{n-j+1} = \Pi_{i = j+1}^n (1 - \frac{1}{2i \log i})^2 \]
As $n \rightarrow \infty$, the above sum clearly goes to 
\begin{equation}
\lambda_2^2 < \pi, \xi_2^2>\int_0^\infty e^{- x} dx. \label{pharaoh}
\end{equation}
Thus, irrespective of $n_0$, the (nonrandom part of the) coefficients of $- \frac{1}{2} t_1^2$ and $- \frac{1}{2} t_2^2$ go to constants as $n \rightarrow \infty$.
At this point note that as we made $n \rightarrow \infty$ the coefficient of $X_{n_0}$ in the characteristic function $t_1 \sqrt{f_{n - n_0 + 1}}$ goes to zero and similarly for the coefficient of $Y_{n_0}$, which is $t_1 \sqrt{g_{n - n_0 + 1}}$. Thus, fixing $n_0$, as we let $n \rightarrow \infty$, the characteristic function does not have $X_{n_0}, Y_{n_0}$ and the nonrandom part of the coefficients of $- \frac{1}{2} t_1^2$ and $- \frac{1}{2} t_2^2$ go to constants independent of $n_0$. This takes care of the sum of $C_{n-j}$'s, $j = n_0, n_0 + 1, \cdots , n$, as we make $n \rightarrow \infty$.
 
{\bf Step 2:} The calculations here will go into the upper bound for the sum of $R_n$'s. The (nonrandom part of the) coefficient of $- (1/2) t_1 t_2$ is 
\begin{equation}
 \sum_{n_0}^n h_{n-j+1} \lambda_1 \lambda_2 \frac{<\pi, \xi_1 \xi_2>}{\sqrt{j+1} \sqrt{(j+1) \log (j+1)}},\label{tia}
\end{equation}
where  \[ h_{n-j+1} = \Pi_{i = j+1}^n (1 - \frac{1}{2i \log i}) (1 - \frac{\frac{1}{2} - \lambda_1}{i}).\] Clearly \[ j_{n-j+1} \leq  \Pi_{i = j+1}^n (1 - \frac{\frac{1}{2} - \lambda_1}{i}),\] and combining the $\sqrt{j+1}$ of $\sqrt{(j+1) \log (j+1)}$ with the other $\sqrt{j+1}$, it is seen that the term (\ref{tia}) is less than 
\[ \frac{1}{\sqrt{\log (n_0 + 1)}} \sum_{n_0}^n  \Pi_{i = j+1}^n (1 - \frac{\frac{1}{2} - \lambda_1}{i}) \lambda_1 \lambda_2 <\pi, \xi_1 \xi_2> . \frac{1}{j+1},\]
which goes to 
\begin{equation}
\frac{1}{\sqrt{\log (n_0 + 1)}} \lambda_1 \lambda_2 <\pi, \xi_1 \xi_2> \int_0^\infty e^{-(\frac{1}{2} - \lambda_1)x} dx \label{ram}
\end{equation}
as $n \rightarrow \infty$. Actually here in the expansion of $(1 - 1/(2i \log i)) (1 - ((1/2) - \lambda_1)/i)$ the important contribution comes from $1 - ((1/2) - \lambda_1)/i$, which can later be compared with the comments following theorem 5.1.

The coefficient of $- (1/2) t_1 t_3$ is (we approximate $\Pi_0^j (1 + \lambda_3/(l+1)) \sim j^{\lambda_3}$), \[ \sum_{n_0}^n f_{n-j+1} \lambda_1 \lambda_3 \frac{<\pi, \xi_1 \xi_3 >}{\sqrt{j+1} j^{\lambda_3}},\] where \[ f_{n-j+1} = \Pi_{i = j+1}^n (1 - \frac{\frac{1}{2} - \lambda_1}{i}). \] Following the argument of the previous paragraph, as we let $n \rightarrow \infty$ this coefficient is less than 
\begin{equation}
\frac{1}{n_o^{\lambda_3 - 1/2}} \lambda_1 \lambda_3 <\pi, \xi_1 \xi_3>  \int_0^\infty e^{-(\frac{1}{2} - \lambda_1)x} dx.\label{rame}
\end{equation}
Similarly as $n \rightarrow \infty$, the (nonrandom part of the) coefficient of $- (1/2) t_2 t_3$ is less than 
\begin{equation}
\frac{\sqrt{(n_0 + 1) \log (n_0 + 1)}}{n_o^{\lambda_3}}\lambda_2 \lambda_3 <\pi, \xi_2 \xi_3> \int_0^\infty e^{-x/2} dx.\label{ramesis}
\end{equation} 
Also note that when we iterate backwards the coefficient of $Z_{n_0}$ is still $t_3$ and keeping $n_0$ fixed as we let $n \rightarrow \infty$ the (nonrandom part of the) coefficient of $- \frac{1}{2} t_3^2$ goes to \begin{equation} \sum_{n_0}^\infty \lambda_3^2 \frac{<\pi, \xi_3^2 >}{(j+1)^{2\lambda_3}}. \label{kiki} \end{equation} 
Thus, fixing $n_0$, the sum of $-t_1 t_2, - t_1 t_3, - t_2 t_3, - \frac{1}{2} t_3^2$, multiplied by their respective (constant part of the) coefficients, is bounded by a constant $F_{n_0}$ as we let $n \rightarrow \infty$. The exact form of $F_{n_0}$ is easily obtained from (\ref{ram}), (\ref{rame}), (\ref{ramesis}) and (\ref{kiki}), however for us the important observation will be that $F_{n_0} \rightarrow 0$ as we later make $n_0 \rightarrow \infty$. 

{\bf Step 3:} The calculations here will go into the upper bound for the sum of $R_n$'s. We now concentrate on the random terms. First note that \[ \sup_{n_0 \leq n < \infty} || \frac{\mathbf{W}_n^\prime}{n+1} - \pi ||, \] where $||.||$ denotes the maximum, is a bounded random variable that converges to $0$ a.s. Also $X_n/\sqrt{n} = \mathbf{W}_n^\prime \xi_1/n$ is bounded by a constant and converges to $0$ a.s. as $n_0 \rightarrow \infty$, hence the same holds for \[ \sup_{n_0 \leq n < \infty} X_n^2/n.\] These two observations show that when we iterate backwards the random terms in the coefficient of $- t_1^2/2$ contribute a random variable less in absolute value than
\begin{equation}
const. \Big{\{} \sup_{n_0 \leq n < \infty} || \frac{\mathbf{W}_n^\prime}{n+1} - \pi || + \sup_{n_0 \leq n < \infty} X_n^2/n \Big{\}} \sum_{n_0}^n f_{n-j+1}  \frac{1}{j+1}.\label{abar}
\end{equation}
The $const.$ term here is an upper bound for $e^{C_{n- n_0} + \cdots + C_n}$ and all the terms in step 2 are also to be multiplied by this. Recall that fixing $n_0$ as we make $n \rightarrow \infty$, the sum  $\sum_{n_0}^n f_{n-j+1} \frac{1}{j+1}$ converges to an integral (see (\ref{pharao}), so that the above sum is bounded by a constant for all $n$), showing that 
as we make $n \rightarrow \infty$ keeping $n_0$ fixed, the contribution of the random terms to the coefficient of $- t_1^2/2$ is bounded by a bounded random variable. This random variable is a constant times the conditional expectation of the random term in (\ref{abar}) given ${\cal F}_{n_0}$, which as we can see converges almost surely to $0$ as we later make $n_0 \rightarrow \infty$.

Similarly, for the other terms involving $Y_n$ and $Z_n$, we use that $\sqrt{\log n / n} Y_n$ and $Z_n/n^{1 - \lambda_3}$ are bounded random variables. Then exactly as in the previous paragraph and following the calculations
 leading to (\ref{pharao}), and the other coefficients (\ref{pharaoh}), (\ref{ram}), (\ref{rame}) and (\ref{ramesis}) we see that fixing $n_0$ as we let $n \rightarrow \infty$, the contribution of the random terms is bounded by a bounded random variable, say the conditional expectation given ${\cal F}_{n_0}$ of a certain $G_{n_0}$ (which goes to $0$ almost surely as we later make $n_0 \rightarrow \infty$). 

\subsection{Completion of proof} 

 Let us now write $H_{n_0} = F_{n_0} + G_{n_0}$, that is the remainder term is bounded by the sum of a constant and a random term uniformly in $n$. Notice that $H_{n_0}$ is actually ${\cal F}_\infty$ measurable and in the calculations what we really use is its conditional expectation given ${\cal F}_{n_0}$. Combining steps 1,2 and 3, fixing $n_0$ as we make $n \rightarrow \infty$, we get from (\ref{chess}) and the previous subsection
\begin{eqnarray}
&& \limsup_{n \rightarrow \infty} |E\{e^{i (t_1 X_{n} + t_2 Y_{n} + t_3 Z_{n})}|{\cal F}_{n_0} \} - 
 e^{i t_3 Z_{n_0}}  e^{- \frac{\sigma_1^2}{2} t_1^2 - \frac{\sigma_2^2}{2} t_2^2}| \nonumber \\
 &\leq& E \{ H_{n_0}|{\cal F}_{n_0} \} + const. \sum_{n_0}^\infty \frac{1}{j^{3/2}},
\end{eqnarray}
with $\sigma_1^2$ and $\sigma_2^2$ coming from (\ref{pharao}) and (\ref{pharaoh}) respectively. Taking expectation and using $ |E V| =  |E E \{ V |{\cal F}_{n_0} \}| \leq E |E \{  V |{\cal F}_{n_0}\}| $, for any integrable random variable $V$,  we get
\begin{equation}
\limsup_{n \rightarrow \infty} |E e^{i (t_1 X_{n} + t_2 Y_{n} + t_3 Z_{n})} - E e^{i t_3 Z_{n_0}}  e^{- \frac{\sigma_1^2}{2} t_1^2 - \frac{\sigma_2^2}{2} t_2^2}| 
 \leq E H_{n_0} + const. \sum_{n_0}^\infty \frac{1}{j^{3/2}}.
\end{equation}
Now $Z_n$ is a martingale, and in the appendix we show that $Z_n$ is $L^2$-bounded, so that $Z_n$ converges to some $Z$ a.s. In the calculation so far $n_0$ is arbitrary. We now let $n_0 \rightarrow \infty$, recalling that the nonrandom $F_{n_0}$ converges to $0$ and that the bounded random variable $G_{n_0}$ also converges to $0$ almost surely from step 3, to get the limiting characteristic function \[ E e^{i t_3 Z} e^{- \frac{\sigma_1^2}{2} t_1^2 - \frac{\sigma_2^2}{2} t_2^2}.\] This shows that $Z$ is independent of $X, Y$, and that $X$ and $Y$ are independent normals. \hfill $\Box$  

\section{Case of real vectors}
\label{sec.real}

In the previous sections we have considered linear combinations corresponding to eigenvectors. To consider general vectors we need the Jordan form of the irreducible replacement matrix. For simplicity we assume that there are only three real eigenvalues. 
However now there exists a nonsingular matrix $\mathbf{T}$ such that 
\[ \mathbf{T}^{-1} \mathbf{R} \mathbf{T} = \left( \begin{array}{cccc} 
1& & & \\ & \Lambda_1 & & \\ & & \Lambda_2 & \\ & & & \Lambda_3 \end{array} \right),\]
where \[ \Lambda_i = \left( \begin{array}{cccc} \lambda_i & 1 & 0 & \\ 0 & \lambda_i & 1 & \\ & & \ddots & \\ & & & \lambda_i \end{array} \right).\] 
Let us consider the case of $\Lambda_1$. The dimension be $d_1$. Then the vectors $\xi_1 = (1, 0, 0, ...)^\prime, \xi_2 = (0, 1, 0, ...)^\prime$, $\cdots$, $\xi_{d_1} = (0, 0, ... , 1)^\prime$ transform according to the equations $\Lambda_1 \xi_1 = \lambda_1 \xi_1, \Lambda_1 \xi_2 = \xi_1 + \lambda_1 \xi_2, \Lambda_1 \xi_3 = \xi_2 + \lambda_1 \xi_3, \cdots$, i.e. in matrix form $\Lambda_1 (\xi_1, \xi_2 , \cdots , \xi_{d_1}) =   (\xi_1, \xi_2 , \cdots , \xi_{d_1}) \Lambda_1$. Denoting the matrix of $\xi_i$'s for the three matrices $\Lambda_1, \Lambda_2, \Lambda_3$ by $\Xi_1, \Xi_2, \Xi_3$ respectively (and necessarily adding $0$'s for the other components) we have \[ \left( \begin{array}{cccc} 
1& & & \\ & \Lambda_1 & & \\ & & \Lambda_2 & \\ & & & \Lambda_3 \end{array} \right) (\mathbf{u}: \Xi_1 : \Xi_2 : \Xi_3) =  (\mathbf{u}: \Xi_1 : \Xi_2 : \Xi_3) \left( \begin{array}{cccc} 
1& & & \\ & \Lambda_1 & & \\ & & \Lambda_2 & \\ & & & \Lambda_3 \end{array} \right),\]
where $\mathbf{u}$ denotes the vector $(1, 0, \cdots)$ of dimension $1 + d_1 + d_2 + d_3$. It may be noticed that $(\mathbf{u}: \Xi_1 : \Xi_2 : \Xi_3)$ is the identity matrix written in a suitable form.

In reality we have to work with not the above matrix of $\Lambda_i$'s, but the stochastic matrix $\mathbf{R}$. In that case, using the above mentioned Jordan decomposition of $\mathbf{R}$, we have to use the vectors 
$\mathbf{T} (\mathbf{u} : \Xi_1 : \Xi_2 : \Xi_3)$, and the equation 
\[ \mathbf{R} \mathbf{T} (\mathbf{u} : \Xi_1 : \Xi_2 : \Xi_3) = \mathbf{T} (\mathbf{u}: \Xi_1 : \Xi_2 : \Xi_3)\left( \begin{array}{cccc} 
1& & & \\ & \Lambda_1 & & \\ & & \Lambda_2 & \\ & & & \Lambda_3 \end{array} \right) .\] As $\mathbf{R}$ has principal eigenvalue 1 corresponding to the eigenvector $\mathbf{1}$ consisting of $1$'s, we have $\mathbf{T} \mathbf{u} = \mathbf{1}$. This implies a trivial limit for $\mathbf{W}_n^\prime \mathbf{T} \mathbf{u}/(n+1)$. However the limits for the other linear combinations corresponding to $\mathbf{W}_n^\prime \mathbf{T} \Xi_i, i = 1, 2, 3$, are nontrivial and are discussed in the next three subsections. For simplicity with a slight abuse of notation we shall use the same notation $\Xi_i$ to denote $\mathbf{T} \Xi_i$.
   
Notice that we can write $\Lambda_i = \lambda_i I_i + F_i$ where $F_i$ is a nilpotent matrix. The presence of this nilpotent $F_i$ changes our calculations in the previous section at certain places and we'll discuss how. We first note that $\mathbf{W}_{n+1}^\prime \Xi_i = \mathbf{W}_{n}^\prime \Xi_i + \chi_{n+1}^\prime \mathbf{R} \Xi_i = \mathbf{W}_{n}^\prime \Xi_i + \chi_{n+1}^\prime \Xi_i \Lambda_i$ (remember the abuse of notation mentioned before). We give the most important contributions, the higher order terms have been ignored for notational simplicity.
  
\subsection{$\lambda_1 < 1/2$}

For notational simplicity from now on we shall restrict ourselves to the highest order terms significant for the results to hold, and this will be denoted by the notation $\sim$. For $\lambda < 1/2$, the approximation $\sqrt{n/(n+1)} \sim (1 - 1/(2n))$ gives
\begin{equation}
 E \{ \frac{\mathbf{W}_{n+1}^\prime \Xi_1}{\sqrt{n+1}}|{\cal F}_n \} \sim \frac{\mathbf{W}_{n}^\prime \Xi_1}{\sqrt{n}} (I_1 - \frac{\frac{1}{2} I_1 - \Lambda_1}{n}),\label{nita}
\end{equation}
leading to the product terms when iterating backwards. On the other hand the approximate form leading to the explicit computations for the conditional characteristic function comes from 
\begin{eqnarray}
&&  \frac{\mathbf{W}_{n+1}^\prime \Xi_1}{\sqrt{n+1}} - E \{ \frac{\mathbf{W}_{n+1}^\prime \Xi_1}{\sqrt{n+1}}|{\cal F}_n \} \nonumber  \\
&\sim& \frac{1}{\sqrt{n+1}} (\chi_{n+1}^\prime - \frac{\mathbf{W}_{n+1}^\prime}{n+1}) \Xi_1 \Lambda_1. \label{nab}
\end{eqnarray}
As before the most inportant contribution in the conditional covariance comes from the first term of the above. Notice that $E \{ \chi_{n+1} \chi_{n+1}^\prime|{\cal F}_n\}$ consists only of diagonal terms and is thus approximately (using Gouet's strong law) $D_{\pi}$, meaning the diagonal matrix with components of $\pi$, namely $\pi_1, \pi_2, \cdots$, as diagonals. This gives for the conditional covariance of (\ref{nab}) the approximate expression \[ \frac{1}{n+1} \Lambda_1^\prime \Xi_1^\prime D_\pi \Xi_1 \Lambda_1 .\]   This when iterated backwards with terms coming from (\ref{nita}), leads to the limiting covariance matrix of the asymptotically normal $\mathbf{W}_n^\prime \Xi_1/\sqrt{n}$, given by 
\begin{eqnarray}
&&\lim_{n\rightarrow \infty} \sum_{n_0}^n \frac{1}{j+1} \Pi_{i = j+1}^n (I_1 - \frac{\frac{1}{2} I_1 - \Lambda_1}{i})^\prime \Lambda_1^\prime \Xi_1^\prime D_\pi \Xi_1 \Lambda_1 \Pi_{i = j+1}^n (I_1 - \frac{\frac{1}{2} I_1 - \Lambda_1}{i}) \nonumber \\
&& = \int_0^\infty e^{- (\frac{1}{2} I_1 - \Lambda_1)^\prime s} \Lambda_1^\prime \Xi_1^\prime D_\pi \Xi_1 \Lambda_1  e^{- (\frac{1}{2} I_1 - \Lambda_1) s} ds,
\end{eqnarray} 
which can be compared with (\ref{pharao}) for the case of eigenvector $\xi_1$.

\subsection{$\lambda_2 = 1/2$}

In this case the norming for the central limit theorem is $\sqrt{n \log^{2 d_2 - 1} n}$ where $d_2$ is the dimension of $\Lambda_2$. The reason for the $2 d_2 - 1$ power will be clear towards the end. First note the approximation \[ \sqrt{ \frac{n \log^{2 d_2 - 1} n}{(n+1) \log^{2 d_2 - 1} (n+1)}} \sim (1 - \frac{1}{2n}) ( 1 - \frac{2d_2 - 1}{2 n \log n}).\] With this we get 
\begin{eqnarray}
&& E \{ \frac{\mathbf{W}_{n+1}^\prime \Xi_2}{\sqrt{(n+1) \log^{2 d_2 - 1} (n+1)}}  |{\cal F}_n \} \nonumber \\ 
&\sim& \frac{\mathbf{W}_n^\prime \Xi_2}{ \sqrt{n \log^{2 d_2 - 1} n} }(1 - \frac{1}{2n}) ( 1 - \frac{2d_2 - 1}{2 n \log n}) \nonumber \\
&+& \frac{\mathbf{W}_n^\prime}{n+1} \frac{\Xi_2 \Lambda_2}{\sqrt{n \log^{2 d_2 - 1} n} } \nonumber \\
&=& \frac{\mathbf{W}_n^\prime \Xi_2}{ \sqrt{n \log^{2 d_2 - 1} n} } ( I_2(1 - \frac{2 d_2 - 1}{2 n \log n}) + \frac{F_2}{n}), \label{ranjan.1}
\end{eqnarray}
where we have crucially used the form of $\Lambda_2$ to cancel the $1/(2n)$'s occuring with opposite signs. This $F_2$ plays an important role in the computations later explaining the $2 d_2 - 1$ power. On the other hand the martingale terms for the covariance computations come from 
\begin{eqnarray}
&&  \frac{\mathbf{W}_{n+1}^\prime \Xi_2}{\sqrt{(n+1) \log^{2 d_2 - 1} (n+1)}} - E \{ \frac{\mathbf{W}_{n+1}^\prime \Xi_2}{\sqrt{(n+1) \log^{2 d_2 - 1} (n+1) }}|{\cal F}_n \} \nonumber  \\
&\sim& \frac{1}{\sqrt{(n+1)\log^{2 d_2 - 1} (n+1) }} (\chi_{n+1}^\prime - \frac{\mathbf{W}_{n+1}^\prime}{n+1}) \Xi_2 \Lambda_2. \label{abhishek.1}
\end{eqnarray}
This gives for the conditional covariance of (\ref{abhishek.1}) the approximate expression \[ \frac{1}{(n+1) \log^{2 d_2 - 1} (n+1) } \Lambda_2^\prime \Xi_2^\prime D_\pi \Xi_2 \Lambda_2 .\]   This when iterated backwards with terms coming from (\ref{ranjan.1}), leads to the limiting covariance matrix of the asymptotically normal $\mathbf{W}_n^\prime \Xi_1/\sqrt{n \log^{2 d_2 - 1} n}$, given by 
\begin{eqnarray}
&\lim_{n \rightarrow \infty}& \sum_{n_0}^n \frac{1}{(j+1) \log^{2 d_2 - 1} (j+1)} \Pi_{i = j+1}^n ( I_2(1 - \frac{2 d_2 - 1}{2 i \log i}) + \frac{F_2}{i})^\prime \nonumber \\ && \Lambda_2^\prime \Xi_2^\prime D_\pi \Xi_2 \Lambda_2    \Pi_{i = j+1}^n ( I_2(1 - \frac{2 d_2 - 1}{2 i \log i}) + \frac{F_2}{i}).  \label{mlk}
\end{eqnarray}
$F_2$ being nilpotent, in the above products only a few terms will be nonzero. The consideration of the limits of the nonzero terms will explain the $\log^{2 d_2 - 1} n$ term in the norming. We illustrate the case of $d_2 = 2$ first. The general case will follow similarly. 

In this case $F_2^2 = \mathbf{0}$. Thus the terms are of three types, (i) $F_2$ is omitted from both sides (ii) $F_2$ occurs on the left but not on the right and conversely (iii)$F_2$ occurs on both sides. Keeping the matrices in order, we bring the constant terms together for the purpose of taking limit.

(i) In this case the products from both sides together reduces to \[ \exp \{ - (2 d_2 - 1) \sum_{j+1}^n \frac{1}{i \log i}\} \sim \exp \{ -(2 d_2 - 1) (\log \log n - \log \log (j+1)) \}.\] Combining this with the $1/(j+1) \log^{2 d_2 - 1} (j+1)$ outside leads to cancellation of  $\log^{2 d_2 - 1} (j+1)$, leaving us with the sum \[ \sum_{n_0}^n \frac{1}{j+1} \frac{1}{\log^{2 d_2 - 1} n} \sim \frac{\log n - \log n_0}{\log^{2 d_2 - 1} n},\] which goes to $0$ as $n \rightarrow \infty$.

(ii) $F_2$ occurs on the left but not on the right. The index from which $F_2$ is taken be $k_1$ which lies between $j+1$ and $n$. Then the constant terms from the products reduce to 
\begin{eqnarray}
&& \exp \{ - \frac{2 d_2 - 1}{2} \sum_{j+1}^{k_1} \frac{1}{i \log i}\} \frac{1}{k_1} \exp \{ - \frac{2 d_2 - 1}{2} \sum_{k_1 + 1}^{n} \frac{1}{l \log l}\} \nonumber \\
&\times & \exp \{ - \frac{2 d_2 - 1}{2} \sum_{j+1}^n \frac{1}{i \log i}\} \nonumber \\
&\sim& \frac{1}{k_1} \exp \{ -(2 d_2 - 1) (\log \log n - \log \log (j+1)) \}. \nonumber 
\end{eqnarray}
Hence combining with the terms outside the sum reduces to 
\begin{eqnarray}
&& \frac{1}{\log^{2 d_2 - 1} n} \sum_{n_0}^n \frac{1}{j+1} \sum_{j+1}^n \frac{1}{k_1}  \nonumber \\
&\sim&  \frac{1}{\log^{2 d_2 - 1} n} \sum_{n_0}^n \frac{1}{j+1} (\log n - \log (j+1)) \nonumber \\
&\sim& \frac{1}{\log^{2 d_2 - 1} n} \{ \log n (\log n - \log n_0) - \frac{1}{2}(\log^2 n - 
\log^2 n_0) \}, \nonumber
\end{eqnarray}
which goes to $0$ as $n \rightarrow \infty$.

(iii) In this case we need one $k_1$ from the left and one $k_2$ from the right corresponding to the indices of $F_2$'s from the left and the right. From the previous calculations the final computation reduces to  
\begin{eqnarray}
&& \frac{1}{\log^{2 d_2 - 1} n} \sum_{n_0}^n \frac{1}{j+1} \sum_{j+1}^n \frac{1}{k_1}  \sum_{j+1}^n \frac{1}{k_2}  \nonumber \\
&\sim& \frac{1}{\log^{2 d_2 - 1} n} \sum_{n_0}^n \frac{1}{j+1} (\log n - \log (j+1))^2 \nonumber \\
&\rightarrow& \frac{1}{3} \nonumber 
\end{eqnarray}
as $n \rightarrow \infty$.

Thus in the case $d_2 = 2$, only the terms corresponding to the highest power of $F_2$ from both sides survives. Exactly the same thing happens for general $d_2$. For example when $d_2 = 3$ the following product corresponds to the highest power of $F_2$ namely $F_2^2$: \[ \Pi_{j+1}^{k_1} (1 - \frac{2 d_2 - 1}{2 i \log i}) \frac{1}{k_1 + 1} 
\Pi_{k_1 +1}^{k_2} (1 - \frac{2 d_2 - 1}{2 l \log l}) \frac{1}{k_2 + 1} \Pi_{k_2 + 1}^{n} (1 - \frac{2 d_2 - 1}{2 m \log m}),\] and the analysis proceeds as before.
The limiting covariance matrix can obtained from (\ref{mlk}), however our main focus is the independence issue. 
Similar calculations can be done using exponentiation technique which is applied
in Section \ref{sec.comp}.

\subsection{$\lambda_3 > 1/2$} 

We expect to get an $L^2$-bounded martingale sequence. Notice first that $E \{ \mathbf{W}_{n+1}^\prime \Xi_3|{\cal F}_n \} = \mathbf{W}_n^\prime \Xi_3 (I_3 + \frac{1}{n+1} \Lambda_3)$. Hence the martingale sequence we work with is 
\begin{equation}
\mathbf{Z}_n =  \mathbf{W}_n^\prime \Xi_3 \{ \Pi_0^{n - 1}( I_3 + \frac{1}{j+1} \Lambda_3)\}^{-1} = \mathbf{W}_n^\prime \Xi_3 \mathbf{A}_n^{-1}. \label{pen}
\end{equation}
The following calculation is similar to the calculation in the appendix and we have used some approximations for notational convenience. $\mathbf{Z}_n$ satisfies the following equation
\begin{eqnarray}
\mathbf{Z}_{n+1} - \mathbf{Z}_n &=& \mathbf{W}_n^\prime \Xi_3 ((I_3 + \frac{1}{n+1} \Lambda_3)^{-1} - I_3) \mathbf{A}_n^{-1} + \chi_{n+1}^\prime \Xi_3 \Lambda_3 \mathbf{A}_{n+1}^{-1} \nonumber \\
&\sim&  - \frac{1}{n+1} \mathbf{Z}_n \mathbf{A}_n \Lambda_3 \mathbf{A}_n^{-1} + \chi_{n+1}^\prime \Xi_3 \mathbf{A}_n^{-1} \mathbf{A}_n \Lambda_3 \mathbf{A}_{n}^{-1} \nonumber \\
&\sim& - \frac{1}{n+1} \mathbf{Z}_n \Lambda_3 + \chi_{n+1}^\prime \Xi_3 \mathbf{A}_n^{-1} \Lambda_3,
\end{eqnarray} 
noting that $\mathbf{A}_n, \mathbf{A}_n^{-1}$ and $\Lambda_3$ commute.
To prove $L^2$-boundedness consider $E E \{ \mathbf{Z}_{n+1} \mathbf{Z}_{n+1}^\prime|{\cal F}_n \}$. Using the martingale property and the above decomposition it follows that
\begin{eqnarray}
E \{ \mathbf{Z}_{n+1} \mathbf{Z}_{n+1}^\prime|{\cal F}_n \} &\sim& \mathbf{Z}_n \mathbf{Z}_n^\prime - \frac{1}{(n+1)^2} \mathbf{Z}_n \Lambda_3 \Lambda_3^\prime \mathbf{Z}_n^\prime \nonumber \\
&& + E \{ \chi_{n+1}^\prime \Xi_3 \mathbf{A}_n^{-1}  \Lambda_3 \Lambda_3^\prime (\mathbf{A}_n^{-1})^\prime \Xi_3^\prime \chi_{n+1} |{\cal F}_n\} \nonumber \\
&\leq& \mathbf{Z}_n \mathbf{Z}_n^\prime (1 - \frac{\beta}{(n+1)^2}) \nonumber \\
&& + Tr \{ \Xi_3 \mathbf{A}_n^{-1}  \Lambda_3 \Lambda_3^\prime (\mathbf{A}_n^{-1})^\prime \Xi_3^\prime E \{ \chi_{n+1} \chi_{n+1}^\prime|{\cal F}_n\} \},
\end{eqnarray}
where $\beta$ denotes the minimum eigenvalue of $\Lambda_3 \Lambda_3^\prime$ and we have used properties of the trace of a matrix. Approximating $E \{ E \{ \chi_{n+1} \chi_{n+1}^\prime|{\cal F}_n\} \}$ by $D_\pi$, further expectation of the above inequality gives 
\begin{eqnarray}
 E \mathbf{Z}_{n+1} \mathbf{Z}_{n+1}^\prime &\leq& E \mathbf{Z}_n \mathbf{Z}_n^\prime (1 - \frac{\beta}{(n+1)^2}) \nonumber \\
 && + const. Tr \{ \Xi_3 \mathbf{A}_n^{-1}  \Lambda_3 \Lambda_3^\prime (\mathbf{A}_n^{-1})^\prime \Xi_3^\prime D_\pi \}.\label{haiti}
\end{eqnarray}
We need to find the order of the last matrix so that the above equation can be iterated as in the one dimensional case of the appendix, giving $L^2$-boundedness of $\mathbf{Z}_n$. Essentially we need to show that the terms of $\mathbf{A}_n^{-1}$ are $O(n^{- \lambda_3} \log^{d_3 - 1} n)$. We do this for $d_3 = 3$, the calculations can be generalised. First we get a formula for $\mathbf{A}_n$. In the following the exponential approximation to $(1 + 1/(j+1))$ has been used for even small values of $j$ for notational simplicity, the exact analysis is similar. First we have
\begin{eqnarray}
\mathbf{A}_n &=& \Pi_1^n (I_3 + \frac{1}{j+1} \Lambda_3) \nonumber \\
&=& \Pi_1^n ( I_3( 1 + \frac{\lambda_3}{j+1}) + \frac{1}{j+1} F_3). \label{mci}
\end{eqnarray}
Using commutativity, the fact that $F_3^2 = \mathbf{0}$, and \[ \Pi_1^n (1 + \frac{\lambda_3}{j+1}) \sim e^{\lambda_3 \log n} \sim n^{\lambda_3},\] $\mathbf{A}_n$ can be written as \[ \mathbf{A}_n = n^{\lambda_3} (I_3 + c_1(n) F_3 + c_2(n) F_3^2).\] $c_1(n)$ comes from 
\begin{eqnarray}
n^{\lambda_3} c_1(n) &=& \sum_{k_1 = 1}^n \Pi_1^{k_1} (1 + \frac{\lambda_3}{j+1}) \frac{1}{k_1 + 1} \Pi_{k_1 + 1}^{n} (1 + \frac{\lambda_3}{l+1}) \nonumber \\
&\sim& \sum_{k_1 = 1}^n e^{\lambda_3 \log k_1} \frac{1}{k_1 + 1} e^{\lambda_3 ( \log n - \log k_1 )} \nonumber \\
&\sim& e^{\lambda_3 \log n }  \sum_{k_1 = 1}^n \frac{1}{k_1 + 1} \nonumber \\
&\sim& n^{\lambda_3} \log n. 
\end{eqnarray}  
A similar computation yields 
\begin{eqnarray}
n^{\lambda_3} c_2(n) &=& \sum_{k_1 = 1}^n \sum_{k_2 = k_1 + 1}^{n} \Pi_1^{k_1} (1 + \frac{\lambda_3}{j+1}) \frac{1}{k_1 + 1} \Pi_{k_1 + 1}^{k_2} (1 + \frac{\lambda_3}{l+1}) \nonumber \\
&& \times \frac{1}{k_2 + 1} \Pi_{k_2 + 1}^n (1 + \frac{\lambda_3}{m+1}) \nonumber \\
&\sim& n^{\lambda_3} \sum_{k_1 = 1}^n \frac{1}{k_1 + 1} \sum_{k_2 = k_1 + 1}^{n} \frac{1}{k_2 + 1} \nonumber \\
&\sim& n^{\lambda_3} \sum_{k_1 = 1}^n \frac{1}{k_1 + 1} (\log n - \log k_1) \nonumber \\
&\sim& n^{\lambda_3} \log^2 n.
\end{eqnarray}
From the above two estimations one can now get the following form of \[ \mathbf{A}_n^{-1}
= n^{-\lambda_3} (I_3 - c_1(n) F_2 + d_2(n) F_2^2),\] where $d_2(n)$ is easily obtained by multiplying the formulas for $\mathbf{A}_n$ and $\mathbf{A}_n^{-1}$ and demanding that the product be $I_3$. This shows that the terms of $\mathbf{A}_n^{-1}(\mathbf{A}_n^{-1})^\prime$ are $O(n^{-2\lambda_3} \log^4 n)$. With this one can go back to (\ref{haiti}) to prove $L^2$-boundedness. Similar calculations can be done using exponentiation technique which is applied
in Section \ref{sec.comp}.

Then the analysis of section 2 proceeds to show independence of the weak limits (strong limit for $\mathbf{Z}_n$). We may state the analogue of theorem 1.1 as follows:
\begin{theorem}
In the case all eigenvalues are real, we consider the linear combinations corresponding to the vectors as identified at the beginning of section 3. The weak limits of the normalized linear combinations corresponding to eigenvalues $\lambda < 1/2, \lambda = 1/2$ and $\lambda > 1/2$ are independent.    
\end{theorem} 
For the different eigenvalues all of which are less than 1/2, there may be dependence among the weak limits coming from the Jordan blocks for different eigenvalues (see theorem 5.1 later). For real $\lambda = 1/2$ there is only one Jordan block (the situation for complex $\lambda$ with real part 1/2 is somewhat different). For $\lambda > 1/2$ the weak limits coming from the Jordan blocks corresponding to different $\lambda$'s are correlated, one instance of this has been computed in the appendix.
Thus they are not independent in general.  

\section{Complex eigenvalues} 
\label{sec.comp}

For complex eigenvalues we consider another canonical form which is similar to the Jordan canonical form. This form comes from considering the real vectors coming from the real and imaginary parts of the complex vectors corresponding to the complex Jordan form. Special cases of this decomposition has been studied in Smythe \cite{smythe}. 
We first consider three types of
eigenvalues, one of each type as before (i.e. with real part less than 1/2. equal to 1/2, and greater than 1/2). 
There exists a nonsingular matrix $\mathbf{S}$ such that 
\[ \mathbf{S}^{-1} \mathbf{R} \mathbf{S} = \left( \begin{array}{cccc} 
1& & & \\ & \Lambda_{c1} & & \\ & & \Lambda_{c2} & \\ & & & \Lambda_{c3} \end{array} \right),\]
where \[ \Lambda_{ci} = \left( \begin{array}{cccc} B_i & I &  & \\   & B_i & \ddots & \\ & & \ddots & I \\ & & & B_i \end{array} \right)\] 
and 
 \[ B_i = \left( \begin{array}{cccc} \lambda_{ir} & \lambda_{ic} \\ - \lambda_{ic} & \lambda_{ir} \end{array}  \right),\] 
$I$ is a 2-dimensional identity matrix 
and 
rest of the elements are 0.
Let the dimension of $\Lambda_{ci}$ be $2 d_{ci}$. 
As before we partition the matrix $\mathbf{S} \mathbf{I}$ (this $\mathbf{I}$ has dimension $1 + 2(d_{c1} + d_{c2} + d_{c3}$)) into a vector of ones and $S_i$, $i=1, 2, 3$ with number of columns in $S_i$ equal to $2 d_{ci}$. These vectors give us the linear combinations.

Notice that, here we can write 
\[ \Lambda_{ci} = \lambda_{ir} I_{ci} + \lambda_{ic} C_{ci} + F_{ci}, \] 
where $I_{ci}$ is an identity matrix of dimension $2d_{ci}$,
 $C_{ci}$ is a block diagonal matrix of same dimension as $\Lambda_{ci}$.
Each block, say $D_i$, is of dimension 2, where
 \[ D_i = \left( \begin{array}{cccc} 0 & 1 \\ - 1 & 0 \end{array} \right).\] 
$F_{ci}$ is a nilpotent matrix of order $d_{ci}$, i.e.,
$F_{ci}^{d_{ci}}= {\bf 0}$ and $d_{ci}$ is the least such integer.

First observe that the rotation matrix $D_i$ satisfies $D_i^2 = - I, D_i^3 = - D_i, D_i^4 = I, \cdots$,
where $I$ is the identity matrix of same dimension as $D_i$.
Also, it is to be noted that the matrices $I_{ci}$, $C_{ci}$ and $F_{ci}$
commute with each other.
Thus,
\begin{eqnarray}
\label{exp.1}
e^{k_1 I_{ci} + k_2 C_{ci} + k_3 F_{ci}} 
&=& e^{k_1 I_{ci}} e^{k_2 C_{ci}} e^{k_3 F_{ci}} \nn\\
&=& e^{k_1} I_{ci}[\cos(k_2) I_{ci} + \sin(k_2) C_{ci}] [\sum_{j=1}^{d_{ci}-1} (k_3^j F_{ci}^j)/j!] .
\end{eqnarray}
We will mention briefly how the proof of theorem \ref{th.1} go for the
complex roots with the presence of the nilpotent matrix and the rotation matrix. 
We note that
$\mathbf{W}_{n+1}^\prime S_i = \mathbf{W}_{n}^\prime S_i + \chi_{n+1}^\prime \mathbf{R} S_i = \mathbf{W}_{n}^\prime S_i + \chi_{n+1}^\prime S_i \Lambda_{ci}$. We give the most important contributions, the higher order terms have been ignored for notational simplicity.
  
\subsection{$\lambda_{1r} < 1/2$}

In this case, since $\sqrt{n/(n+1)} \sim (1 - 1/(2n))$,
 it is to be noted that,
\begin{equation}
 E \{ \frac{\mathbf{W}_{n+1}^\prime S_1}{\sqrt{n+1}}|{\cal F}_n \} \sim \frac{\mathbf{W}_{n}^\prime S_1}{\sqrt{n}} (I_{c1} - \frac{\frac{1}{2} I_{c1} - \Lambda_{c1}}{n}) .\label{nita.1}
 \end{equation}
Now iterating backwards we get the product terms as before. 
Thus,
\begin{eqnarray}
&&  \frac{\mathbf{W}_{n+1}^\prime S_1}{\sqrt{n+1}} - E \{ \frac{\mathbf{W}_{n+1}^\prime S_1}{\sqrt{n+1}}|{\cal F}_n \} \nonumber  \\
&\sim& \frac{1}{\sqrt{n+1}} (\chi_{n+1}^\prime - \frac{\mathbf{W}_{n}^\prime}{n+1}) S_1 \Lambda_{c1}. \label{nab.1}
\end{eqnarray}
As before the most important contribution in the conditional covariance comes from the first term of the above. Notice that $E \{ \chi_{n+1} \chi_{n+1}^\prime|{\cal F}_n\}$ consists only of diagonal terms and is thus approximately (using Gouet's strong law) $D_{\pi}$. This gives for the conditional covariance of (\ref{nab.1}) the approximate expression \[ \frac{1}{n+1} \Lambda_{c1}^\prime S_1^\prime D_\pi S_1 \Lambda_{c1} .\]   This when iterated backwards with terms coming from (\ref{nita.1}), leads to the limiting covariance matrix of the asymptotically normal $\mathbf{W}_n^\prime S_1/\sqrt{n}$, given by 
\begin{eqnarray}
&\lim_{n\rightarrow \infty}& \sum_{n_0}^n \frac{1}{j+1}\Big{\{} \Pi_{i = j+1}^n (I_{c1} - \frac{\frac{1}{2} I_{c1} - \Lambda_{c1}}{i})^\prime \times \nonumber \\
&&  \Lambda_{c1}^\prime S_1^\prime D_\pi S_1 \Lambda_{c1} \Pi_{i = j+1}^n (I_{c1} - \frac{\frac{1}{2} I_{c1} - \Lambda_{c1}}{i}) \Big{\}} \nonumber \\
&=& \int_0^\infty e^{- (\frac{1}{2} I_{c1} - \Lambda_{c1})^\prime s} \Lambda_{c1}^\prime S_1^\prime D_\pi S_1 \Lambda_{c1}  e^{- (\frac{1}{2} I_{c1} - \Lambda_{c1}) s} ds,\label{lim.1}
\end{eqnarray} 
which can be compared with (\ref{pharao}) for the case of eigenvector $\xi_1$.
From the calculation in (\ref{exp.1}), it can be seen that
$$
e^{- (\frac{1}{2} I_{c1} - \Lambda_{c1}) s}
=e^{- (\frac{1}{2}-\lambda_{1r})s} I_{c1} [\cos(s\lambda_{1c}) I_{c1} + \sin(s\lambda_{1c})C_{c1}] [\sum_{j=1}^{d_{c1}-1} (s F_{c1})^j/j!] ,
$$
which is an integrable function, and hence (\ref{lim.1}) is finite.

\subsection{$\lambda_{2r} = 1/2$}

In this case the norming for the central limit theorem is $\sqrt{n \log^{2 d_{c2}-1} n}$ where $2d_{c2}$ is the dimension of $\Lambda_{2c}$. From the calculation of the  covariance matrix the reason for the  $2 d_{c2}-1$ power of the 
the logarithm will be clear. The approximation \[ \sqrt{ \frac{n \log^{2 d_{c2}-1} n}{(n+1) \log^{2 d_{c2}-1} (n+1)}} \sim (1 - \frac{1}{2n}) ( 1 - \frac{2d_{c2}-1}{2 n \log n})\] leads to 
\begin{eqnarray}
&& E \{ \frac{\mathbf{W}_{n+1}^\prime S_2}{\sqrt{(n+1) \log^{2 d_{c2}-1} (n+1)}}  |{\cal F}_n \} \nonumber \\ 
&\sim& \frac{\mathbf{W}_n^\prime S_2}{ \sqrt{n \log^{2 d_{c2}-1} n} }(1 - \frac{1}{2n}) ( 1 - \frac{2d_{c2}-1}{2 n \log n}) \nonumber \\
&+& \frac{\mathbf{W}_n^\prime}{n+1} \frac{S_2 \Lambda_2}{\sqrt{n \log^{2 d_{c2}-1} n} } \nonumber \\
&=& \frac{\mathbf{W}_n^\prime S_2}{ \sqrt{n \log^{2 d_{c2}-1} n} } ( I_{c2}(1 - \frac{2 d_{c2}-1}{2 n \log n}) + \frac{\lambda_{2c}}{n} C_{c2} + \frac{F_{c2}}{n}), \label{ranjan.2}
\end{eqnarray}
where the form of $\Lambda_{2c}$ is used to cancel the $1/(2n)$'s occuring with opposite signs. We later discuss the role of $C_{c2}$ and $F_{c2}$ in the computations that explains the power of the logarithm. 
Notice that the computation of the covariance matrix depends on the martingale terms
\begin{eqnarray}
&&  \frac{\mathbf{W}_{n+1}^\prime S_2}{\sqrt{(n+1) \log^{2 d_{c2}-1} (n+1)}} - E \{ \frac{\mathbf{W}_{n+1}^\prime S_2}{\sqrt{(n+1) \log^{2 d_{c2}-1} (n+1) }}|{\cal F}_n \} \nonumber  \\
&\sim& \frac{1}{\sqrt{(n+1)\log^{2 d_{c2}-1} (n+1) }} (\chi_{n+1}^\prime - \frac{\mathbf{W}_{n+1}^\prime}{n+1}) S_2 \Lambda_{c2}. \label{abhishek.2}
\end{eqnarray}
Thus, the approximate expression
for the conditional covariance of (\ref{abhishek.2}) is found as
 \[ \frac{1}{(n+1) \log^{2 d_{c2}-1} (n+1) } \Lambda_{c2}^\prime S_2^\prime D_\pi S_2 \Lambda_{c2} .\]   Iterating backwards with terms coming from (\ref{ranjan.2}) leads to the limiting covariance matrix of the asymptotically normal $\mathbf{W}_n^\prime S_2/\sqrt{n \log^{2 d_{c2}-2} n}$, given by 
\begin{eqnarray}
&\lim_{n \rightarrow \infty}& \sum_{n_0}^n \frac{1}{(j+1) \log^{2 d_{c2}-1} (j+1)} \nonumber \\
&& \Big{\{} \Pi_{i = j+1}^n ( I_{c2}(1 - \frac{2 d_{c2}-1}{2 i \log i}) +
\frac{\lambda_{2c}}{i} C_{c2}^{\prime} + \frac{F_{c2}}{i})^\prime \Lambda_{c2}^\prime S_2^\prime D_\pi S_2 \Lambda_{c2} \nonumber \\ 
&& \times \Pi_{i = j+1}^n ( I_{c2}(1 - \frac{2 d_{c2}-1}{2 i \log i}) + \frac{\lambda_{2c}}{i} C_{c2} + \frac{F_{c2}}{i}) \Big{\}}.  \label{mlk.1}
\end{eqnarray}
We shall now use exponentiation to simplify the calculations (the same could be done in section 3.2). Observing that,
\begin{eqnarray}
&&
\Pi_{i = j+1}^n ( I_{c2}(1 - \frac{2 d_{c2}-1}{2 i \log i}) + \frac{\lambda_{2c}}{i} C_{c2} + \frac{F_{c2}}{i}) \nn\\
& \sim &
\Pi_{i = j+1}^n e^ {- \frac{2 d_{c2}-1}{2 i \log i}I_{c2} + \frac{\lambda_{2c}}{i} C_{c2} + \frac{F_{c2}}{i}} \nn\\
& = &
 e^ {- \sum_{i = j+1}^n \frac{2 d_{c2}-1}{2 i \log i}I_{c2} + \sum_{i = j+1}^n\frac{\lambda_{2c}}{i} C_{c2} + \sum_{i = j+1}^n\frac{F_{c2}}{i}} \nn\\
& \sim &
 e^ {- \frac{2 d_{c2}-1}{2} \log \frac{\log n}{\log (j+1)}I_{c2} +
  \lambda_{2c} C_{c2} \log (\frac{n}{j+1})  + F_{c2} \log (\frac{n}{j+1})} \nn\\
& = &
 e^ {- \frac{2 d_{c2}-1}{2} \log \frac{\log n}{\log (j+1)}}I_{c2} 
 \ e^{C_{c2} \lambda_{2c} \log (\frac{n}{j+1})} \ e^{F_{c2} \log (\frac{n}{j+1})} \nn\\
& = &
 e^ {- \frac{2 d_{c2}-1}{2} \log \frac{\log n}{\log (j+1)}}I_{c2} 
 \ [I_{c2} \cos(\lambda_{2c} \log (\frac{n}{j+1})) + C_{c2} \sin(\lambda_{2c} \log (\frac{n}{j+1}))] \nn\\
&&   [ \sum_{k=1}^{d_{c2} - 1} (F_{c2} \log (\frac{n}{j+1}))^k/k!] .    
\label{prdt}
\end{eqnarray}
Combining the contribution of the term in (\ref{prdt}) to the two sides of 
 (\ref{mlk.1}) we get
\begin{eqnarray}
&& \frac{1}{(j+1) \log^{2 d_{c2}-1} (j+1)} \Pi_{i = j+1}^n ( I_{c2}(1 - \frac{2 d_{c2}-1}{2 i \log i}) +
\frac{\lambda_{2c}}{i} C_{c2} + \frac{F_{c2}}{i})^\prime \nonumber \\ 
&& \Lambda_{c2}^\prime S_2^\prime D_\pi S_2 \Lambda_{c2}  \Pi_{i = j+1}^n ( I_{c2}(1 - \frac{2 d_{c2}-1}{2 i \log i}) + \frac{\lambda_{2c}}{i} C_{c2} + \frac{F_{c2}}{i}) \nn\\
& \sim &
\frac{1}{(j+1) \log^{2 d_{c2}-1} (j+1)} e^ {- (2 d_{c2}-1) \log \frac{\log n}{\log (j+1)}} \nn\\
&& \times [I_{c2} \cos(\lambda_{2c} \log (\frac{n}{j+1})) + C_{c2}^{\prime} \sin(\lambda_{2c} \log (\frac{n}{j+1}))] \nonumber \\
&& \times [ \sum_{k=1}^{d_{c2}-1} (F_{c2}^{\prime} \log (\frac{n}{j+1}))^k/k!] 
\Lambda_{c2}^\prime S_2^\prime D_\pi S_2 \Lambda_{c2} \nonumber \\
&& \times [I_{c2} \cos(\lambda_{2c} \log (\frac{n}{j+1})) + C_{c2} \sin(\lambda_{2c} \log (\frac{n}{j+1}))] \nn\\
&& \times [ \sum_{k=1}^{d_{c2}-1} (F_{c2} \log (\frac{n}{j+1}))^k/k!] \nn
\end{eqnarray}
\begin{eqnarray}
&=&
\frac{1}{(j+1) \log^{2 d_{c2}-1} n}
\ [I_{c2} \cos(\lambda_{2c} \log (\frac{n}{j+1})) + C_{c2}^{\prime} \sin(\lambda_{2c} \log (\frac{n}{j+1}))] \nn\\
&& \ [ \sum_{k=1}^{d_{c2}-1} (F_{c2}^{\prime} \log (\frac{n}{j+1}))^k/k!] \nn\\
&&
\Lambda_{c2}^\prime S_2^\prime D_\pi S_2 \Lambda_{c2} 
\ [I_{c2} \cos(\lambda_{2c} \log (\frac{n}{j+1})) + C_{c2} \sin(\lambda_{2c} \log (\frac{n}{j+1}))] \nn\\
&& \ [ \sum_{k=1}^{d_{c2}-1} (F_{c2} \log (\frac{n}{j+1}))^k/k!] 
\end{eqnarray}
Now arguing as in section \ref{sec.real}
and observing that terms involving $sine$ and $cosine$ are all bounded,
one finds that except the coefficient of the highest power term of $F_{c2}$ i.e. $F_{c2}^{d_{c2} - 1}$, the coefficients of other terms go
to zero when $n \to \infty$. Observe that highest power terms of $F_{c2}$ would be multiplied by 
$\cos^2(\lambda_{2c} \log (\frac{n}{j+1}))$ (i.e. $(1/2) [1 + \cos(2\lambda_{2c} \log (\frac{n}{j+1}))])$, 
$\sin^2(\lambda_{2c} \log (\frac{n}{j+1}))$ i.e. $(1/2) [1 - \cos(2\lambda_{2c} \log (\frac{n}{j+1}))])$,
\ or, terms such as, \
$\sin(\lambda_{2c} \log (\frac{n}{j+1}))\ \cos(\lambda_{2c} \log (\frac{n}{j+1}))$ (i.e. $(1/2) [\sin(2\lambda_{2c} \log (\frac{n}{j+1}))])$,
separately.
Thus, 
the highest power terms of $F_{c2}$ with $sine$ function give the coefficient
\begin{eqnarray}
&& \frac{1}{\log^{2 d_{c2}-1} n} \sum_{j=n_0}^n \frac{1}{j+1} \Big{\{} \frac{(\log n - \log (j+1))^{2d_{c2}-2}}{((d_{c2}-1)!)^2} \nonumber \\ && \times \frac{\sin (2\lambda_{2c}(\log n - \log (j+1)))}{2} \Big{\}} \nonumber \\
&\sim & \frac{1}{\log^{2 d_{c2}-1} n} \int_{0}^{\log n - \log n_0} \frac{u^{2d_{c2}-2}}{((d_{c2}-1)!)^2} \frac{\sin (2\lambda_{2c} u)}{2} \nonumber \\ 
&=& O\left(\frac{(\log n - \log n_0)^{2d_{c2}-2}}{\log^{2 d_{c2}-1} n}\right)
\ \ \rightarrow 0 
\end{eqnarray}
(seen by integration by parts) as $n \rightarrow \infty$.
Similarly, with $cosine$ function, it gives
\begin{eqnarray}
&& \frac{1}{\log^{2 d_{c2}-1} n} \sum_{j=n_0}^n \frac{1}{j+1}\Big{\{} \frac{(\log n - \log (j+1))^{2d_{c2}-2}}{((d_{c2}-1)!)^2} \nonumber \\ 
&& \times \frac{\cos (2\lambda_{2c}(\log n - \log (j+1)))}{2} \Big{\}} \nonumber \\
&\sim & \frac{1}{\log^{2 d_{c2}-1} n} \int_{0}^{\log n - \log n_0} \frac{u^{2d_{c2}-2}}{((d_{c2}-1)!)^2} \frac{\cos (2\lambda_{2c} u)}{2} \nonumber \\
&=& O\left(\frac{(\log n - \log n_0)^{2d_{c2}-2}}{\log^{2 d_{c2}-1} n}\right)
\ \ \rightarrow 0 
\ \ \ \mbox{as} \ \ n \to \infty . 
\end{eqnarray}
Now, the terms that involve multiplying by the \ $1/2$ \ only, give
\begin{eqnarray}
&& \frac{1}{\log^{2 d_{c2}-1} n} \sum_{j=n_0}^n \frac{1}{j+1} \frac{(\log n - \log (j+1))^{2d_{c2}-2}}{((d_{c2}-1)!)^2} \frac{1}{2} \nonumber \\
&\sim & \frac{1}{\log^{2 d_{c2}-1} n} \int_{0}^{\log n - \log n_0} \frac{u^{2d_{c2}-2}}{((d_{c2}-1)!)^2} \frac{1}{2} \nonumber \\
&=& \frac{1}{2 (2d_{c2}-1) \ ((d_{c2}-1)!)^2}
 \left(\frac{(\log n - \log n_0)^{2d_{c2}-1}}{\log^{2 d_{c2}-1} n}\right) \nn\\
&\rightarrow& \frac{1}{2 (2d_{c2}-1) \ ((d_{c2}-1)!)^2} \label{rie}
\end{eqnarray}
as $n \rightarrow \infty$. Thus the limiting covariance matrix can obtained from (\ref{mlk.1}).
Notice that it does not involve $\lambda_{2c}$.

\subsection{$\lambda_{3r} > 1/2$} 

Here also, we show that $\mathbf{W}_n^\prime S_3 A_n^{-1}$ is
an $L^2$-bounded martingale sequence, where
$A_n = \Pi_0^{n - 1}( I_{c3} + \frac{1}{j+1} \Lambda_{c3})$.
 Notice first that $E \{ \mathbf{W}_{n+1}^\prime S_3|{\cal F}_n \} = \mathbf{W}_n^\prime S_3 (I_{c3} + \frac{1}{n+1} \Lambda_{c3})$. Hence the martingale 
$\mathbf{Z}_n$ satisfies the following equation
\begin{eqnarray}
\mathbf{Z}_{n+1} - \mathbf{Z}_n &=& \mathbf{W}_n^\prime S_3 ((I_{c3} + \frac{1}{n+1} \Lambda_{c3})^{-1} - I_{c3}) \mathbf{A}_n^{-1} + \chi_{n+1}^\prime S_3 \Lambda_{c3} \mathbf{A}_{n+1}^{-1} \nonumber \\
&\sim&  - \frac{1}{n+1} \mathbf{Z}_n \mathbf{A}_n \Lambda_{c3} \mathbf{A}_n^{-1} + \chi_{n+1}^\prime S_3 \mathbf{A}_n^{-1} \mathbf{A}_n \Lambda_{c3} \mathbf{A}_{n}^{-1} \nonumber \\
&\sim& - \frac{1}{n+1} \mathbf{Z}_n \Lambda_{c3} + \chi_{n+1}^\prime S_3 \mathbf{A}_n^{-1} \Lambda_{c3},
\end{eqnarray} 
since $\mathbf{A}_n, \mathbf{A}_n^{-1}$ and $\Lambda_{c3}$ commute.
To prove $L^2$-boundedness, first observe
\begin{eqnarray}
E \{ \mathbf{Z}_{n+1} \mathbf{Z}_{n+1}^\prime|{\cal F}_n \} &\sim& \mathbf{Z}_n \mathbf{Z}_n^\prime - \frac{1}{(n+1)^2} \mathbf{Z}_n \Lambda_{c3} \Lambda_{c3}^\prime \mathbf{Z}_n^\prime \nonumber \\
&& + E \{ \chi_{n+1}^\prime S_3 \mathbf{A}_n^{-1}  \Lambda_{c3} \Lambda_{c3}^\prime (\mathbf{A}_n^{-1})^\prime S_3^\prime \chi_{n+1} |{\cal F}_n\} \nonumber \\
&\leq& \mathbf{Z}_n \mathbf{Z}_n^\prime (1 - \frac{\beta_{c3}}{(n+1)^2}) \nonumber \\
&& + Tr \{ S_3 \mathbf{A}_n^{-1}  \Lambda_{c3} \Lambda_{c3}^\prime (\mathbf{A}_n^{-1})^\prime S_3^\prime E \{ \chi_{n+1} \chi_{n+1}^\prime|{\cal F}_n\} \},
\end{eqnarray}
where $\beta_{c3}$ denotes the minimum eigenvalue of $\Lambda_{c3} \Lambda_{c3}^\prime$. 
Approximating \newline $E \{ E \{ \chi_{n+1} \chi_{n+1}^\prime|{\cal F}_n\} \}$ by $D_\pi$, further expectation of the above inequality gives 
\begin{eqnarray}
 E \mathbf{Z}_{n+1} \mathbf{Z}_{n+1}^\prime &\leq& E \mathbf{Z}_n \mathbf{Z}_n^\prime (1 - \frac{\beta_{c3}}{(n+1)^2}) \nonumber \\
 && + const. Tr \{ S_3 \mathbf{A}_n^{-1}  \Lambda_{c3} \Lambda_{c3}^\prime (\mathbf{A}_n^{-1})^\prime S_3^\prime D_\pi \}.\label{haiti.1}
\end{eqnarray}
We now find the order of the last matrix so that the above equation can be iterated as in the one dimensional case of the appendix, giving $L^2$-boundedness of $\mathbf{Z}_n$. 
We show this by showing that
the terms of $\mathbf{A}_n^{-1}$ are $O(n^{- \lambda_3} \log^{d_3 - 1} n)$. 
\begin{eqnarray}
\mathbf{A}_n &=& \Pi_1^n (I_{c3} + \frac{1}{j} \Lambda_{c3}) \nonumber \\
&=& \Pi_1^n ( I_{c3}( 1 + \frac{\lambda_{3r}}{j}) + \frac{\lambda_{3r}}{j} C_{c3} + \frac{1}{j} F_{c3}). \label{mci.1}
\end{eqnarray}
Using commutativity of $I_{c3}, C_{c3}$ and $F_{c3}$ and the fact that $F_{c3}^{d_{c3}} = \mathbf{0}$, 
$\mathbf{A}_n$ can be approximated as \[ \mathbf{A}_n 
\sim e^{\lambda_{3r} \log n I_{c3} + C_{c3} \log n + F_{c3} \log n}
\]
Hence 
\begin{eqnarray}
\mathbf{A}_n^{-1}
&\sim& e^{-\lambda_{3r} \log n} I_{c3} 
 [\cos(-\lambda_{3r}\log n) I_{c3} + \sin(-\lambda_{3r}\log n) C_{c3}] \nonumber \\
&& \times [ \sum_{k=1}^{d_{c3}-1} F_{c3}^k (-\log n)^k/k!] \nonumber \\
&\sim& n^{\lambda_{3r}} \log^{d_{c3}-1} n. 
\end{eqnarray}  
Thus, $\mathbf{A}_n^{-1} = O(n^{-\lambda_{3r}} \log^{d_{c3}-1} n) $ and from (\ref{haiti.1}) one gets $L^2$-boundedness of $\mathbf{Z}_n$. 

Then the analysis of section 2 proceeds to show independence of the weak limits (with strong limit for $\mathbf{Z}_n$). We may state the analogue of theorem 1.1 as follows:
\begin{theorem}
In the case eigenvalues are complex, we consider the linear combinations corresponding to the vectors as identified at the beginning of section 4. The weak limits of the normalized linear combinations corresponding to eigenvalues $Re(\lambda) < 1/2, Re(\lambda) = 1/2$ and $Re(\lambda) > 1/2$ are independent.    
\end{theorem} 
For the different eigenvalues all of which have real parts less than 1/2, there may be dependence among the weak limits coming from the (modified) Jordan blocks for different eigenvalues (see theorem 5.1 later). For $Re(\lambda) = 1/2$, there may be different (modified) Jordan blocks corresponding to different $Im(\lambda)$. However, inside $Re(\lambda) = 1/2$, the weak limits coming from (modified) Jordan blocks of different dimensions are not independent, in general. 
For $\lambda > 1/2$ the weak limits coming from the Jordan blocks corresponding to different $\lambda$'s are correlated (similar to the real eigenvalue case computed in the appendix).
Thus they are not independent in general.

\section{General case}
\label{sec.gen}

In the general case we decompose the replacement matrix into a (modified) Jordan form 
as in the previous two sections. That is, corresponding to real eigenvalues we take the form as in section 3, and corresponding to complex eigenvalues by considering the real and imaginary parts of vectors we take the form from section 4.
Without loss of generality, now we can consider only the real parts of the eigenvalues,  and the linear combinations will come from the (modified) Jordan form. 

There are now three types of blocks, for $Re(\lambda) < 1/2$,
for $Re(\lambda) = 1/2$ and the last type is for $Re(\lambda) > 1/2$. According to our previous notation, there exists a nonsingular matrix $\mathbf{M}$ such that 
\[ \mathbf{M}^{-1} \mathbf{R} \mathbf{M} = \left( \begin{array}{cccc} 
1& & & \\ & G_{1} & & \\ & & G_{2} & \\ & & & G_{3} \end{array} \right),\]
where \[ G_{i} = \left( \begin{array}{cccc} \Lambda_{i,1} &  &  & \\  & \Lambda_{i,2} &  & \\ & & \ddots & \\ & & & \Lambda_{i,n_i} \end{array} \right)\] 
and $\Lambda_{i,j}$s are either of the form of $\Lambda_{i}$ as in section
(\ref{sec.real}) or $\Lambda_{ci}$ as in section
(\ref{sec.comp}).
Also notice that, for each $i=1,2,3$,  there is a positive integer 
$0 \leq k_i \leq n_i$ such that $\Lambda_{i,1}, \ldots, \Lambda_{i,k_i}$ blocks correspond to real eigenvalues and the rest of the $n_i - k_i$ blocks 
correspond to complex eigenvalues. It can be observed that $k_2 \leq 1$, and it is also assumed that the blocks inside $G_2$ which have the same dimension (i.e. same $d_2$ or $d_{c2}$) are arranged next to one another and put into one subblocks. 

Let us recall  that the linear combinations come from the columns of $\mathbf{M} \mathbf{I}$ which we write with an abuse of notation as $(\mathbf{1}:M_1:M_2:M_3)$. With appropriate normalizations they decompose into the following three classes, independent in the limit.
\begin{theorem}

{\bf 1. Re($\lambda) <$ 1/2:} For the linear combinations corresponding to columns of $M_1$, the normalization is $\sqrt{n}$ and the limit is normal. The covariance is given by (\ref{lim.1}) with $G_1$ replacing $\Lambda_{c1}$ (and $M_1$ replacing $S_1$) and we have to use the decomposition of $G_1$ combining the features of the real and the complex cases. 

{\bf 2. Re($\lambda$) = 1/2:} Recalling the arrangement inside $G_2$, in this case the linear combinations correspond to columns of $M_2$. For the subblock of $G_2$ having dimension $d_2$ or $d_{c2}$ for the original $\Lambda_{2,k}$'s (of the same dimension), the normalization for the corresponding columns of $M_2$ is $\sqrt{n \log^{2d_{c2}  - 1} n }$ (or   $\sqrt{n \log^{2d_{2}  - 1} n}$ as appropriate) and the limit is normal. 
The limits for different subblocks are not independent, in general, and for each subblock the covariance can be found from (\ref{mlk.1}) to (\ref{rie}) by decomposing the subblock of $G_2$ combining the features of the real and the complex cases (and replacing $S_2$ by the column submatrix of $M_2$ corresponding to the subblock of $G_2$).

{\bf 3. Re($\lambda) > 1/2$:}  For the linear combinations corresponding to columns of $M_3$, $\mathbf{W}_n^\prime M_3 A_n^{-1}$ is
an $L^2$-bounded martingale sequence, where
$A_n = \Pi_0^{n - 1}( I_{3} + \frac{1}{j+1} G_{3})$, and $I_3$ is an identity matrix of the same dimension as $G_3$. The covariance among some of the components of the (almost sure) limit may be nonzero, even though rates are different, implying dependence.  
\end{theorem}

To summarize part one and two of the above theorem, observe that $(1/\sqrt{n})$
is the only normailization constant for the part one, i.e., for $W_n'M_1$
and (not necessarily zero) covariances are obtained between different Jordan
blocks in this part.
Whereas, for part two, let us take
$M_2 = [M_{2, 1} : \ldots : M_{2, n_{2}}]$ where $M_{2, j}$'s correspond to the different Jordan subblocks. Then 
$$W_n'M_2 C_{n2} = \left(W_n'M_{2, 1}, W_n'M_{2, 2}, \ldots, W_n'M_{2, n_{2}}\right) 
 \left( \begin{array}{cccc} C_{n2,1} &  &  & \\  & C_{n2,2} &  & \\ & & \ddots & \\ & & & C_{n2,n_i} \end{array} \right) 
$$
is asymptotically normal with covariance matrix given below, where $C_{n2,j}$ is a
 diagonal matrix of dimension $p_{mj}$ with each entry as 
$(1/\sqrt{n \log^{2d_{mj}-1} n})$.
Here $d_{mj}$ equals to $d_{2j}$ if it corresponds to a real eigenvalue (as in Section \ref{sec.real}),
and it is $d_{c2j}$ if it corresponds to a complex case
(as in Section \ref{sec.comp}),
whereas $p_{mj}$ equals to $d_{2j}$ if it corresponds to a real eigenvalue,
and it is $2d_{c2j}$ if it corresponds to a complex case 
as in Section \ref{sec.comp}). This is a case for asymptotic mixed normality.
In this case,
the typical entries of the limiting covariance matrix of $W_n'M_2 C_n$, say 
$V_{2}$, can be seen as follows, as in (\ref{rie}),
$$V_{2}(j,l) = \frac{1}{2 (d_{mj}+d_{ml}-1) \ ((d_{mj}-1)!(d_{ml}-1)!)} \Lambda_{mj} M_{2,j}' D_{\pi} M_{2,l} \Lambda_{ml}, 
$$ where $\Lambda_{mj}$ is the subblock of $G_2$ corresponding to $M_{2,j}$.

\section{appendix}
\label{sec.append}

Suppose $U_n$ and $V_n$ are normalized linear combinations corresponding to eigenvectors $\xi_3, \xi_4$, with eigenvalues $\lambda_3, \lambda_4$, respectively both of which are real and greater than 1/2. We want to show that the limiting covariance is nonzero. 
This technique has been used in the proof of lemma 3.1 of Freedman \cite{freedman}.
$U_n$ and $V_n$ satisfy the following equations:
\begin{eqnarray}
U_{n+1} - U_n &=&  \lambda_3 \frac{\chi_{n+1}^\prime \xi_3}{\Pi_0^n(1 + \frac{\lambda_3}{j+1})} -\frac{\frac{\lambda_3}{n+1}}{1 + \frac{\lambda_3}{n+1}} U_n, \nonumber \\
V_{n+1} - V_n &=&  \lambda_4 \frac{\chi_{n+1}^\prime \xi_4}{\Pi_0^n(1 + \frac{\lambda_4}{j+1})} -\frac{\frac{\lambda_4}{n+1}}{1 + \frac{\lambda_4}{n+1}} V_n.
\end{eqnarray}
Using the martingale property it follows that
\begin{eqnarray}
E \{ U_{n+1} V_{n+1}|{\cal F}_n \} &=& U_n V_n ( 1 - \frac{\frac{\lambda_3}{n+1}}{1 + \frac{\lambda_3}{n+1}} \frac{\frac{\lambda_4}{n+1}}{1 + \frac{\lambda_4}{n+1}}) \nonumber \\
&& + \frac{\lambda_3 \lambda_4}{\Pi_0^n(1 + \frac{\lambda_3}{j+1})\Pi_0^n(1 + \frac{\lambda_4}{j+1})} < \frac{\mathbf{W}_n}{n+1}, \xi_3 \xi_4>, \nonumber \\
E U_{n+1} V_{n+1} &=& E U_n V_n ( 1 - \frac{\frac{\lambda_3}{n+1}}{1 + \frac{\lambda_3}{n+1}} \frac{\frac{\lambda_4}{n+1}}{1 + \frac{\lambda_4}{n+1}}) \nonumber \\
&& + \frac{\lambda_3 \lambda_4}{\Pi_0^n(1 + \frac{\lambda_3}{j+1})\Pi_0^n(1 + \frac{\lambda_4}{j+1})} < E \frac{\mathbf{W}_n}{n+1}, \xi_3 \xi_4>.
\end{eqnarray}
Notice that by the dominated convergence theorem and Gouet's strong law, $E \frac{\mathbf{W}_n}{n+1}$ converges to $<\pi, \xi_3 \xi_4>$. Iterating the above equation and using that $\Pi_0^n(1 + \frac{\lambda_3}{j+1}) \sim n^{\lambda_3}$  we get (remember $\lambda_3, \lambda_4 > 1/2$) that $E U_n V_n$ converges to a nonzero quantity. In particular the same technique yields the $L^2$-boundedness of $Z_n$ of section 1.

\end{document}